\documentclass[letterpaper, 10 pt, conference]{ieeeconf}  

\IEEEoverridecommandlockouts                              
\usepackage{graphicx}
\usepackage{amsmath} 
\usepackage{mathtools} 
\usepackage{tikz}
\usepackage{pgfplots}
\pgfplotsset{compat=1.18}
\usepackage{amssymb}  

\title{\LARGE \bf
On the threshold of excitable systems: An energy-based perspective
}

\author{Rodolphe Sepulchre and Guanchun Tong
\thanks{This work was supported by the European Research Council under the Advanced ERC Grant Agreement SpikyControl (project number: 101054323).
The work of Guanchun Tong was also supported by Flemish Fund for Scientific Research (FWO) under a PhD Fellowship fundamental research (project number: 1191125N).
}
\thanks{Guanchun Tong is with the Department of Electrical Engineering (ESAT), KU Leuven, Belgium.
        {\tt\small guanchun.tong@esat.kuleuven.be}}%
\thanks{Rodolphe Sepulchre is with the Department of Electrical Engineering(ESAT), KU Leuven, Belgium, and also with the Department of 
Engineering, Control Group, University of Cambridge, UK.
        {\tt\small rodolphe.sepulchre@kuleuven.be}}%
}

\begin{document}

\maketitle
\thispagestyle{empty}
\pagestyle{empty}

\begin{abstract}


A fundamental characteristic of excitable systems is their ability to exhibit distinct subthreshold and suprathreshold behaviors. Precisely quantifying this distinction requires a proper definition of the threshold, which has remained elusive in neurodynamics. In this paper, we introduce a novel, energy-based threshold definition for excitable circuits grounded in dissipativity theory, specifically using the classical concept of required supply. According to our definition, the threshold corresponds to a \textit{local maximum of the required supply}, clearly separating subthreshold passive responses from suprathreshold regenerative spikes. We illustrate and validate the proposed definition through analytical and numerical studies of three canonical systems: a simple RC circuit, the FitzHugh--Nagumo model, and the biophysically detailed Hodgkin--Huxley model.

\end{abstract}

\section{INTRODUCTION}

The division between the discrete and the continuous is pervasive across scientific modeling and technological applications \cite{Lovász2010,SIAMNews2012}. In the context of computing machines, this dichotomy often manifests as the digital versus the analog, each offering unique advantages  \cite{AnalogVersusDigital}. Combining the strengths of both paradigms is central to neuromorphic engineering, wherein spiking behaviors exhibit continuous dynamics yet retain the reliability of discrete events, as these spikes can be individually counted \cite{sepulchre2022spiking}. This mixed property is grounded in neuronal excitability, where each spike represents a discrete transition from continuous, subthreshold activity into a transient, suprathreshold event. Thus, the reliability of spiking systems critically depends on precisely quantifying the boundary between subthreshold continuity and suprathreshold discreteness.

Such dichotomy can only be quantified with a mathematical definition of threshold. But what is an excitability threshold? This question transcends the literature of neurodynamics and has remained elusive. It has a simple answer only for highly idealized models of excitability. For instance, the celebrated leaky integrate-and-fire model \cite{lapicque1907recherches} models the spike as a straightforward reset mechanism in an RC circuit. The voltage is reset whenever it reaches a {\it threshold} value $v_{th}$. However, it has long been acknowledged that this simplistic definition of threshold becomes inadequate in more general models of excitability. The entire first chapter of the standard neurodynamics textbook by Izhikevich \cite{izhikevich2006dynamical} provides an excellent account of the difficulties in defining a ``voltage" threshold or a ``current" threshold for neuronal behaviors. 

The proposed remedy in the same textbook is the general methodology of  neurodynamics: different types of excitable models are classified according to bifurcation analysis. However, this longstanding approach has several limitations. First of all, bifurcation analysis only partially captures the inherent input-output nature of excitability, as it traditionally examines closed dynamical systems using constant inputs as bifurcation parameters. Additionally, excitability intrinsically involves transient phenomena, making steady-state bifurcation analysis indirect and sometimes misleading --- the spiking event of an excitable system is only indirectly related to the steady-state limit cycle that arises from applying a constant current of sufficient amplitude. 

The second limitation concerns scalability. The textbook \cite{izhikevich2006dynamical} classifies no fewer that sixteen different types of bursting models, and the classification heavily relies on phase portraits and time-scale separation. It is unclear how this approach could be generalized to the high-dimensional neuronal models relevant to neuromorphic engineering. 

Finally, and perhaps most importantly, the robustness of a definition of threshold with respect to model uncertainty is a delicate issue.  The mismatch between internal and external robustness of dynamical input-output systems is widely acknowledged as a key question of system theory. It has been primarily studied in the context of equilibrium linear stability, where the mismatch is between the robustness of the eigenspectrum of the internal operator and the sensitivity of the input-output operator. Several earlier studies have illustrated an similar mismatch in neuronal behaviors (see for instance \cite{franci2018robust} for an illustration in the context of single neuron properties). 

As an alternative to the bifurcation theory of excitability threshold, the present paper adopts an energy-based perspective. We aim at characterizing the thresholds in excitable systems described by physical circuits that interconnect storage elements (e.g., capacitors), dissipative elements (e.g., resistive or memristive elements), and active elements such as batteries. The subthreshold behavior corresponds to passive trajectories, which cannot extract energy from internal batteries, therefore they can only dissipate the energy supplied externally. The suprathreshold behavior, in contrast, results from active trajectories capable of extracting energy from internal batteries, thereby enabling a higher storage without additional external supply. We characterize the threshold of such circuits as an energy threshold, by computing the minimal external supply that is required to escape the subthreshold regime of the circuit. This energy-based characterization is firmly rooted in the classical theory of dissipativity \cite{willems1972dissipative}, and leads to a definition of threshold by solving an optimal control problem. The voltage and current thresholds are the corresponding optimal circuit trajectories that bring the system from equilibrium to an event with a minimal supply of energy.

To the best knowledge of the authors, the proposed energy-based perspective is novel. In this preliminary manuscript, we illustrate this definition in the simplest examples of excitable models, and we highlight the  potential of the proposed definition for a general  theory of excitability that would combine desirable physical and algorithmic properties.

The remainder of the article is organized as follows. Section II introduces the Hodgkin–Huxley model as a canonical example of excitable behavior, with a discussion of the difficulties involved in defining a clear threshold. We then outline our energy-based alternative. In Sections III through V, we illustrate the proposed definition using three examples: a linear/cubic RC circuit, the FitzHugh--Nagumo model, and finally a return to the Hodgkin–Huxley model. The FitzHugh--Nagumo model serves as a bridging example between the analytical study of the RC circuit and the numerical results for the Hodgkin–Huxley model.

\section{An energy-based definition of   the threshold in excitable systems}

\subsection{Hodgkin--Huxley model}
The Hodgkin–Huxley (HH) model \cite{Hodgkin1952Quantitative,Keener2009Mathematical} describes neuronal excitability through a nonlinear circuit representation of ion channel dynamics. The membrane potential \( v \) evolves according to:
\begin{equation}
 C \dot{v} = - i_{Na} - i_{K} - i_{L} + i,   
\end{equation}
where \( C \) is the membrane capacitance, \( i \) is the injected current, and \( i_{Na}, i_{K}, i_{L} \) represent sodium, potassium, and leak currents, respectively. The ionic currents are modeled as:
\begin{equation}
   \begin{aligned}
i_{Na} &= \bar{g}_{Na} m^3 h (v - v_{Na}),\\
i_{K} &= \bar{g}_{K} n^4 (v - v_{K}),\\
i_{L} &= \bar{g}_{L} (v - v_{L}),
\end{aligned} 
\end{equation}
with \(\bar{g}_i\) denoting maximal conductances and \(v_i\) reversal potentials. The gating variables \( m, h, n \in [0,1] \) follow first-order kinetics:
\begin{equation}
 \tau_x(v) \dot{x} = -x + x_{\infty}(v), \quad x \in \{m,h,n\},   
 \label{eq:gating variables HH}
\end{equation}
where \(\tau_x(v)\) and \(x_{\infty}(v)\) are time constants and steady-state activations (both nonlinearly dependent on voltage), empirically derived from voltage clamp experiments. Sodium activation (\( m \)) is fast, enabling rapid depolarization, while inactivation (\( h \)) and potassium activation (\( n \)) are slower, ensuring repolarization.

\subsection{Subthreshold, Suprathreshold, and Inhibitory Responses}

\begin{figure}[h]  
     \centering
     \includegraphics[width=0.48\textwidth]{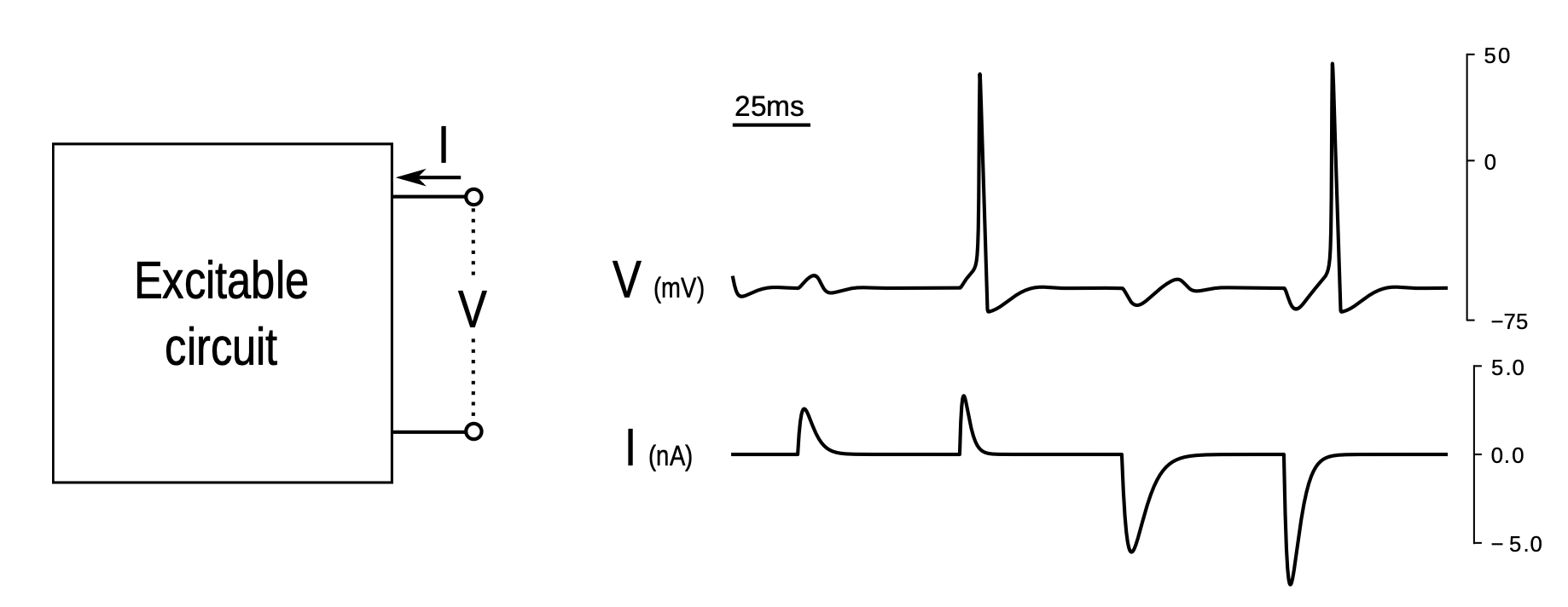} 
     \caption{Subthreshold and suprathreshold behaviors of an excitable circuit. Reproduced from \cite{sepulchre2018excitable} with permission.}
     \label{fig:excitable behavior}
\end{figure}

Excitability is the ability to generate all-or-none spikes when crossing the threshold. Fig.\@~\ref{fig:excitable behavior}  illustrates this behavior: small currents (subthreshold) induce passive voltage changes, while suprathreshold inputs trigger stereotyped spikes.

Notably, hyperpolarizing inputs can also elicit spikes. When inhibitory currents de-inactivate Na\(^{+}\) channels, rebound depolarization is enabled upon its release. Such dual sensitivity underscores that excitability is a dynamic property, not merely a voltage threshold.

Excitability enables neurons to encode information digitally while remaining sensitive to analog input features. In neuromorphic engineering, replicating this property allows energy-efficient event-based information processing across scales \cite{Sepulchre2019Control,ribar2021neuromorphic,Schmetterling2024Neuromorphic}.

\subsection{Why thresholds defy a simple definition?}

\begin{figure}[h]  
     \centering
     \includegraphics[width=0.48\textwidth]{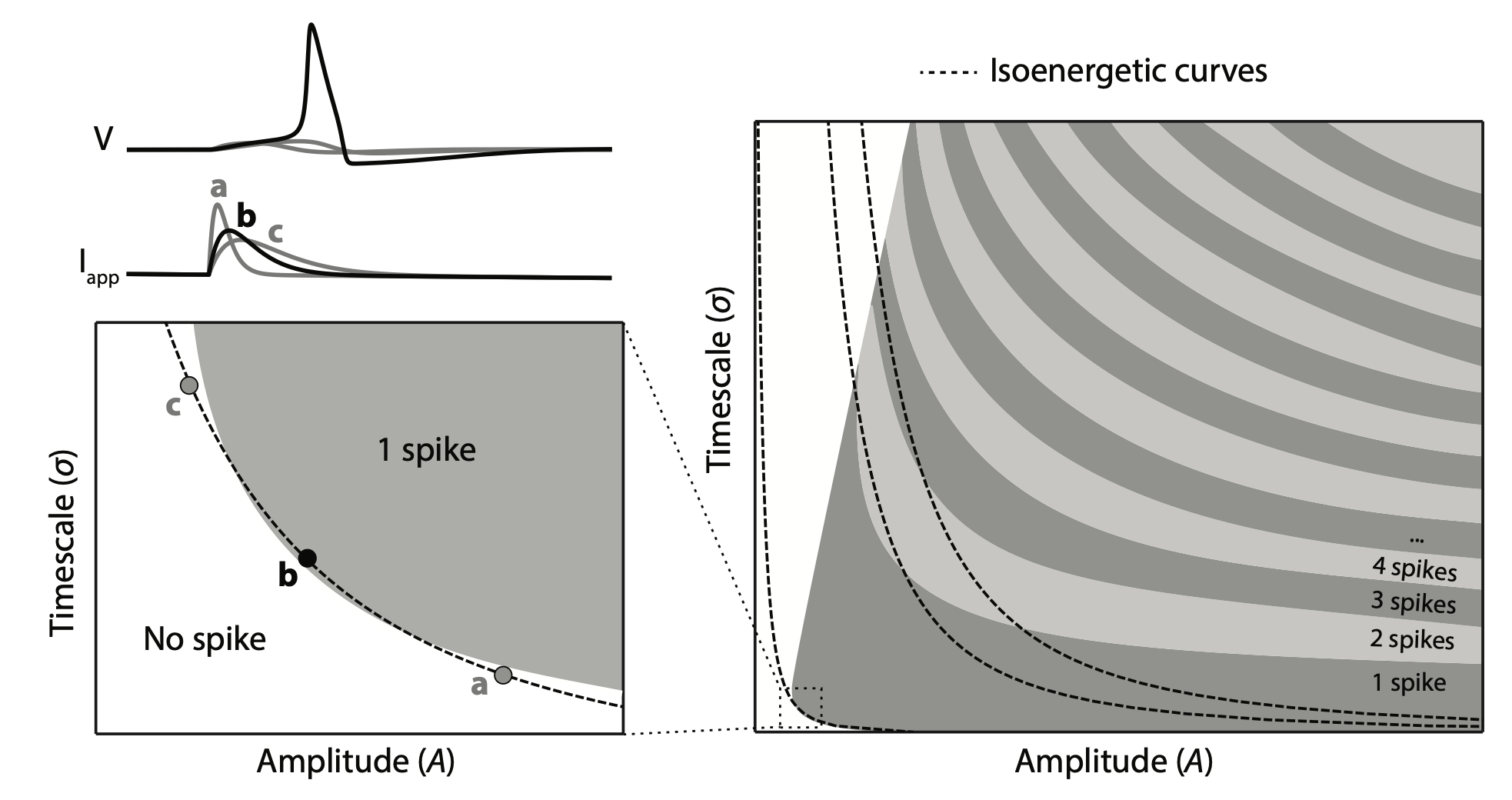} 
     \caption{The subtlety of defining threshold \cite{sepulchre2018excitable}.}
     \label{fig:threshold property}
\end{figure}

Thresholds are often inaccurately perceived as fixed voltage levels. For example, Fig.\@~\ref{fig:threshold property} demonstrates that spike initiation in excitable systems depends intricately on input amplitude  \(A\) and time scale \(\sigma\).  Even on an iso-charge curve (same area under the three currents),  minor variations in temporal profiles significantly affect spiking outcomes. Short-duration inputs fail to activate slow internal gating variables (e.g., \(n\)), while prolonged inputs allow dissipation mechanisms to dominate, preventing spike generation. Therefore, the threshold represents a complex interplay between internal energy dissipation and conservation with external inputs, challenging simplistic definitions.

\subsection{A dissipativity theoretic definition}


Consider the past before crossing the threshold and identify the time reaching the threshold as \(t=0\), the circuit satisfies the energy balance:
\begin{equation}
     \label{eq:energy balance}
\underbrace{\int_{-\infty}^{0}i(t)\,v(t)\,dt}_{\text{external supply}} \;=\; \underbrace{\int_{-\infty}^{0}\sum_{k=1}i_{k}(t)\,v(t)\,dt}_{\text{internal supply + dissipation}} + \underbrace{ S(v(0))}_{\text{storage}},
\end{equation}
where $S(v(0))=\int_{-\infty}^{T} C\dot{v}(t)\, v(t) \, dt = \frac{1}{2}Cv^2(0)$, assuming the capacitor is initially uncharged. Like resistive elements, memristive elements only dissipate energy. Therefore, unlike in the classical framework of dissipativity theory, the complexity of the model lies in characterizing the dissipation rather than the storage, which in the Hodgkin--Huxley model reduces  to the aforementioned capacitive energy.

The required supply is the mimimum energy that must be supplied to the circuit to reach a target state $x^*$ at time $t=0$, starting from the equilibrium state at $t=-\infty$:
\begin{equation}
S_r(x^*)=\inf_{i_{(-\infty, 0]}\colon x(0)=x^*} \int_{-\infty}^{0}i(t)\,v(t)\,dt,
\end{equation}
where \(i_{(-\infty, 0]}\) denotes the input signal to be optimized over. This definition is in fact \(S_{rc}\) in \cite{van_der_schaft_cyclo_dissipativity_2021} with a more symmetrized treatment of required supply and available storage.

We define the energy threshold of an excitable system as a local maximum of the required supply. In other words, this threshold identifies the state that triggers an event without further external energy, allowing the system to reach higher voltages with less external input.


\section{Illustration example: an RC circuit}

The definition of required supply is via solving an optimal control problem where the supply rate $i(t)v(t)$ serves as the running cost. We illustrate this first on a simple RC circuit, highlighting the singular nature of the resulting optimal control problem. This motivates the approximation solution by an ``exponential ansatz''.

\subsection{A singular optimal control problem}

We want to solve the following optimal control problem:
\begin{equation}
  \begin{aligned}
&\text{Minimize} && J \;=\; \int_{-\infty}^{0} i(t)\,v(t)\,dt,\\[4pt]
&\text{subject to} && C\,\dot{v}(t) \;=\; -g(v(t))v(t)\;+\;i(t),\\[4pt]
&&& v(-\infty)\;=\;0,\quad v(0)\;=\;v^*.
\end{aligned}  
\end{equation}

We apply Pontryagin's Minimum Principle by forming the Hamiltonian with costate varialbe $\lambda$:
\begin{equation}
    H(v, i, \lambda) \coloneq iv + \lambda \frac{(-g(v)v+i)}{C}.
\end{equation}


The costate equation is:
\begin{equation}
\label{eq: costate equation}
    \dot{\lambda}(t)= -\frac{\partial H}{\partial v} = \frac{\lambda}{C}(g'(v)v+g(v))-i.
\end{equation}

The necessary condition for optimality is 
\begin{equation}
\label{eq: necessary condition for optimality}
    \frac{\partial H}{\partial i} = v + \frac{\lambda}{C}=0,
\end{equation}
therefore $\lambda(t)=-Cv(t)$.

From  \eqref{eq: costate equation} ,\eqref{eq: necessary condition for optimality} and the state equation, it follows that 
\begin{equation}
    g(v)v+v(g'(v)v+g(v))=0,
\end{equation}
which can be rewritten as 
\begin{equation}
    \frac{d}{dv}[g(v)v^2]=0.
\end{equation}


This is a singular optimal control problem due to linear-in-control cost and linear additive control in the dynamics. To be able to solve a singular optimal control problem, it is typically required to invoke higher-order stationary conditions or leverage more information from the system dynamics. This justifies our approach by using the exponential ansatz for the voltages, which will be discussed in detail in the follwoing.


\subsection{The exponential ansatz}
To illustrate the idea of the exponential ansatz, we now focus on a linear RC circuit. The idea will then be generalized to nonliner bi-stable RC circuits as well as the Fitzhugh--Nagumo model and the Hodgkin--Huxley model.

For a linear RC circuit, the optimal control problem specializes to the following:

\begin{equation}
  \begin{aligned}
&\text{Minimize} && J \;=\; \int_{-\infty}^{0} i(t)\,v(t)\,dt,\\[4pt]
&\text{subject to} && C\,\dot{v}(t) \;=\; -\frac{1}{R}v(t)\;+\;i(t),\\[4pt]
&&& v(-\infty)\;=\;0,\quad v(0)\;=\;v^*.
\end{aligned}  
\end{equation}

Motivated by the system dynamics and the boundary conditions, we make the ansatz about the optimal state trajectory that $v(t)=v^*e^{\alpha t}$; the optimal control trajectory now becomes:
\begin{equation}
    i(t) = C\dot{v}(t) + \frac{1}{R}v(t) = \left(C\alpha + \frac{1}{R}\right)v^* e^{\alpha t}.
\end{equation}

Thus, the supplied energy becomes:
\begin{equation}
    J(v^*, \alpha) = \int_{-\infty}^{0} i(t)\, v(t)\, dt 
= \frac{C\alpha + \frac{1}{R}}{2\alpha} {v^*}^2.
\end{equation}

Therefore, the original infinite-dimensional optimal control problem simplifies to optimizing over one single parameter $\alpha$. In this linear RC case, we have 
\begin{equation}
    S_r(v^*) = \inf_{\alpha} J(v^*, \alpha) = \frac{1}{2} C {v^*}^2.
\end{equation}

The minimum is attained when $\alpha \to \infty$, which corresponds to the optimal policy that applies an impulse current at $t=0$. In this manner, no energy is dissipated. All supplied energy is stored in the capacitor, as illustrated in the optimal value $1/2 C{v^*}^2$. Since the required supplied in this example is convex in $v^*$, there is no local maximum: reaching a higher storage can only be achieved by supplying more energy. This is a energy-based characterization of the fact that the linear RC circuit has no threshold. Its full behavior is {\it subthreshold}, that is, the equilibrium behavior of a dissipative RLC circuit.

\subsection{Generalization to a bi-stable RC circuit}

In the case of nonlinear RC circuit  $C\dot{v}(t) \;=\;-i_{d}(t) + i(t)$, where $i_d(t)=f(v(t))$ has the following
nonlinear characteristic, intersecting the $v-$axis at $0=v_a<v_b<v_c = 1$.

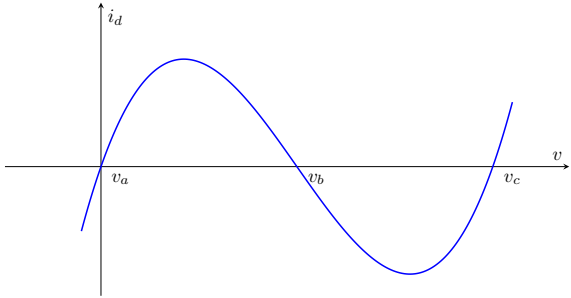
\begin{figure}[h]
\centering
\begin{tikzpicture}[scale=0.72]
 \begin{axis}[
    axis lines=middle,
    xlabel={$v$},
    ylabel={$i_d$},
    samples=200,
    domain=-0.2:4.2,
    ymin=-3, ymax=4,
    xmin=-0.5, xmax=4.3,
    xtick=\empty, 
    ytick={0},
    width=12cm,
    height=7cm,
    enlargelimits=true,
    clip=true,
  ]
    \addplot[thick, blue] {x^3 - 6*x^2 + 8*x};

    \node[below] at (axis cs:0.2,0) {$v_a$};
    \node[below] at (axis cs:2.2,0) {$v_b$};
    \node[below] at (axis cs:4.2,0) {$v_c$};
  \end{axis}
\end{tikzpicture}
\caption{The N-shape i-v curve for the nonlinear resistor. }
\label{fig:bistable RC}
\end{figure}

The required supply reaching voltage $v^*$ at $t=0$ (starting from equilibrium $v(-\infty)=0$) is 
\begin{equation}
     S_r(v^*) = \inf_{i_{(-\infty, 0]}\colon v(0)=v^*} \int_{-\infty}^{0}i(t)\,v(t)\,dt.
\end{equation}

This is equivalent to 
\begin{equation}
     S_r(v^*) \;=\;  \inf_{v_{(-\infty, 0]}\colon v(0)=v^*} \int_{-\infty}^{0} f(v(t))\,v(t)\,dt + \frac{1}{2}C{v^*}^2.
\end{equation}

When $v^* \in [0,v_b]$, the circuit is passive. The minimum required energy is simply $s(v^*)=\frac{1}{2}C{v^*}^2$ by setting the voltage 
level instantly(as in the linear case), therefore the term $\inf \int_{-\infty}^{0} f(v(t))\,v(t)\,dt$ vanishes.

When $v \in (v_b, v_c]$, we enter a region where energy can be extracted internally. One possible optimal control policy consists in setting the voltage level to be 
$v^*$ at the very beginning. The required supply then becomes $-\infty$, which makes $v_b$ a local maximum of the required supply, therefore the  threshold of this bi-stable RC circuit.

The conclusion is that according to the proposed definition, the threshold of a bistable RC circuit is the saddle point of the closed bistable circuit, that is when the input current is zero. The energy threshold is the supply that is required to charge the capacitor up to the initial condition that will trigger an event in the absence of additional supply. For the bistable circuit, the suprathreshold behavior is the switch from the equilibrium of low storage to the equilibrium of high storage.

\section{Thresholds in the  Fitzhugh--Nagumo model}

\subsection{Fitzhugh--Nagumo model}
The FitzHugh--Nagumo (FHN) model \cite{FitzHugh1961Impulses, Nagumo1962Active} serves as a conceptual intermediate between bistable RC circuits and the Hodgkin–Huxley model. The system,

\begin{equation}
   \label{eq:FHN}
  \begin{aligned}
    \epsilon \dot{v} &= -f(v) - w + i,\\
    \dot{w} &= v - \gamma w,
\end{aligned}  
\end{equation}
where 
\begin{equation}
     \quad f(v) = (1 - v)(v - v_b)v.
\end{equation}

features timescale separation ($\epsilon \ll 1$) and a cubic nonlinearity $f(v)$ that induces bistability for $0 < v_b < 1$. In the singular limit ($\epsilon \to 0$), $w$ becomes quasi-static, reducing \eqref{eq:fhn} to a scalar bistable system $\dot{v} = -f(v) - w + i$, analogous to an RC circuit with two stable equilibria (Fig.~ \ref{fig:bistable RC}) . Here, $i$ (with $w$ absorbed into it) acts as a bifurcation parameter, enabling bi-stable/mono-stable transitions. However, unlike scalar RC models, the inclusion of the slow variable  $w$ introduces history-dependent recovery. This allows for the generation of spike events  through the interplay of slow negative feedback and fast positive feedback—a hallmark of neural excitability.

FHN circuit can be regarded as two-dimensional simplification of the HH model's four-dimensional ionic kinetics: the $v$-nullcline's cubic geometry approximates sodium channel activation, while the linear $w$-nullcline mimics slower potassium dynamics. Phase-plane analysis reveals excitable and oscillatory regimes governed by $i$ and $\gamma$, mirroring HH's spiking behaviors without its biophysical complexity. 

\begin{figure}[h]
\centering
\begin{tikzpicture}
\begin{axis}[
    width=0.42\textwidth,
    height=0.42\textwidth,
    xlabel={$v$},
    ylabel={$w$},
    xmin=-0.5, xmax=1.5,
    ymin=-0.5, ymax=0.8,
    grid=major,
    axis equal image
]

\addplot[red, thick, domain=-0.5:1.5, samples=100] 
    {(1 - x)*(x - 0.4)*x}; 
\addplot[blue, thick, domain=-0.5:1.5] 
    {2*x}; 



\node[fill=green, circle, inner sep=2pt] at (axis cs:0,0) {};

\legend{v-nullcline, w-nullcline}
\end{axis}
\end{tikzpicture}
\caption{The red cubic curve represents the $v$-nullcline ($\dot{v}=0$), blue line shows the $w$-nullcline ($\dot{w}=0$). Green circle: stable fixed point }
\end{figure}
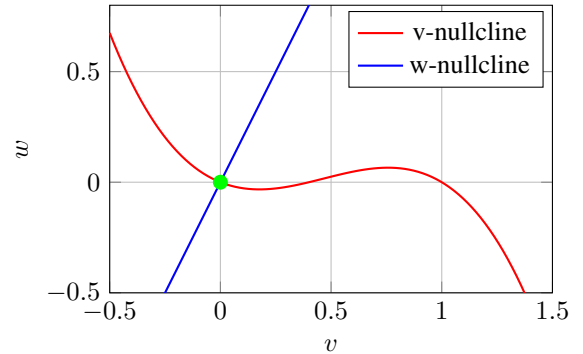
\subsection{Energy thresholds in FHN}
Following the exponential ansatz of the bistabe circuit, we have \(v(t)=Ae^{\alpha t}\). We can solve from \(\dot{w} = v - \gamma w\):
\begin{equation}
    w(t) = \frac{A}{\alpha+\gamma}e^{\alpha t}.
\end{equation}

The applied current is
\begin{equation*}
\begin{aligned}
    i(t) &= \epsilon \dot{v}(t)+f(v(t))+w(t)\\
    &=(\epsilon\alpha-v_b+\frac{1}{\alpha+\gamma}) Ae^{\alpha t} - A^3e^{3\alpha t} + (v_b+1)A^2e^{2\alpha t}.
\end{aligned}
\end{equation*}

The cost functional becomes:
\begin{equation*}
\begin{aligned}
      J(A,\alpha)&=\int_{-\infty}^0 i(t)v(t)dt\\
      &= \frac{1}{2}\epsilon A^2 - \frac{A^4}{4\alpha} + \frac{(v_b+1)A^3}{3\alpha}-\frac{v_bA^2}{2\alpha}+\frac{A^2}{2\alpha(\alpha+\gamma)}.
\end{aligned}
\end{equation*}

\begin{figure}[h]  
     \centering
     \includegraphics[width=0.5\textwidth]{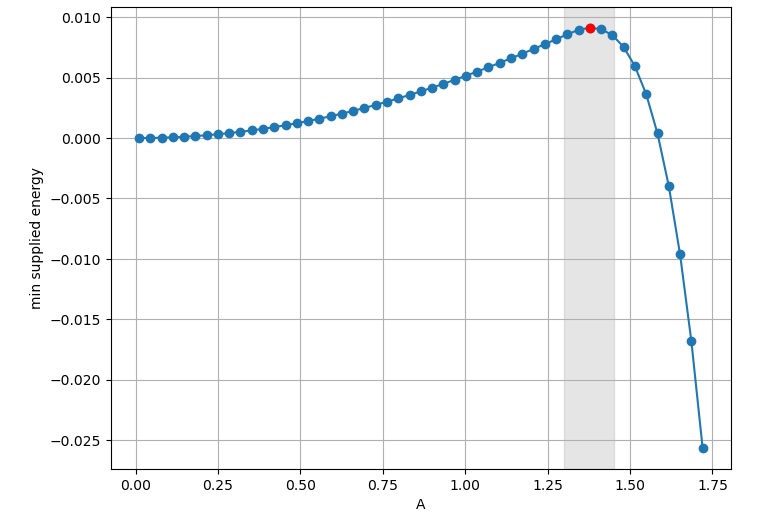} 
     \caption{One threshold of the Fitzhugh--Nagumo circuit.}
     \label{fig:threshold-FHN}
\end{figure}

Fig.~\ref{fig:threshold-FHN} illustrates the minimal energy required as a function of the target voltage \( A \). A local maximum occurs near \(A\approx 1.38\) with \(\alpha \approx 51.10\). Below this threshold voltage, progressively larger energy inputs are necessary to reach higher voltages, reflecting the local passive behavior of the circuit. Once the threshold is exceeded, however, the circuit's negative conductance becomes active, thereby enabling a spiking event without additional energy supply.


\section{Thresholds in the Hodgkin--Huxley Model}

As introduced in Section II, the Hodgkin--Huxley circuit involves four nonlinearly, dynamically coupled state variables: the membrane voltage \( v \) and three gating variables \( m, h\) and \(n\). Solving the full optimal control problem would ordinarily involve all four state variables. However, once a voltage trajectory \(v(t)\) is prescribed, the gating variables become uniquely determined (initial conditions are at equilibrium). Hence, the simplified approach of optimizing over the voltage trajectory alone is justified. It can be regarded as a generalization of the exponential ansatz.

Specifically, we parameterize the voltage trajectory as an exponential function of time: \(v(t) = A e^{\alpha t}, \quad t \leq 0,\)
with \(v(-\infty)=0\) and \(v(0) = A\). For each voltage trajectory defined by parameters \( A \) and \( \alpha \), we update only the gating variables \( m, h, n \) numerically according to \eqref{eq:gating variables HH}, their dynamic equations in the Hodgkin--Huxley model.

\subsection{Numerical procedure: dynamic ``voltage clamp"}
The numerical procedure for computing the minimum supplied energy is outlined as follows:

\begin{itemize}
    \item Choose discretized ranges of \(A\) and \(\alpha\).
    \item For each pair \((A, \alpha)\):
    \begin{enumerate}
        \item Initialize membrane voltage and gating variables at equilibrium values \(v(-\infty), m_{\infty}(v(-\infty))\), \(h_{\infty}(v(-\infty))\), \(n_{\infty}(v(-\infty))\).
        \item For each time step ( from \(t = -\infty\) to \(t = 0\)):
        \begin{enumerate}
            \item Force \( v(t) = A e^{\alpha t} \).
            \item Update gating variables \(m, h, n\) numerically.
            \item Calculate ionic currents and dissipated energy.
        \end{enumerate}
        \item Compute total supplied energy as the sum of dissipated ionic energy and stored capacitive energy.
    \end{enumerate}
    \item For each \(A\), identify the minimum energy over all  \(\alpha\).
\end{itemize}

\begin{figure}[h]  
     \centering
     \includegraphics[width=0.5\textwidth]{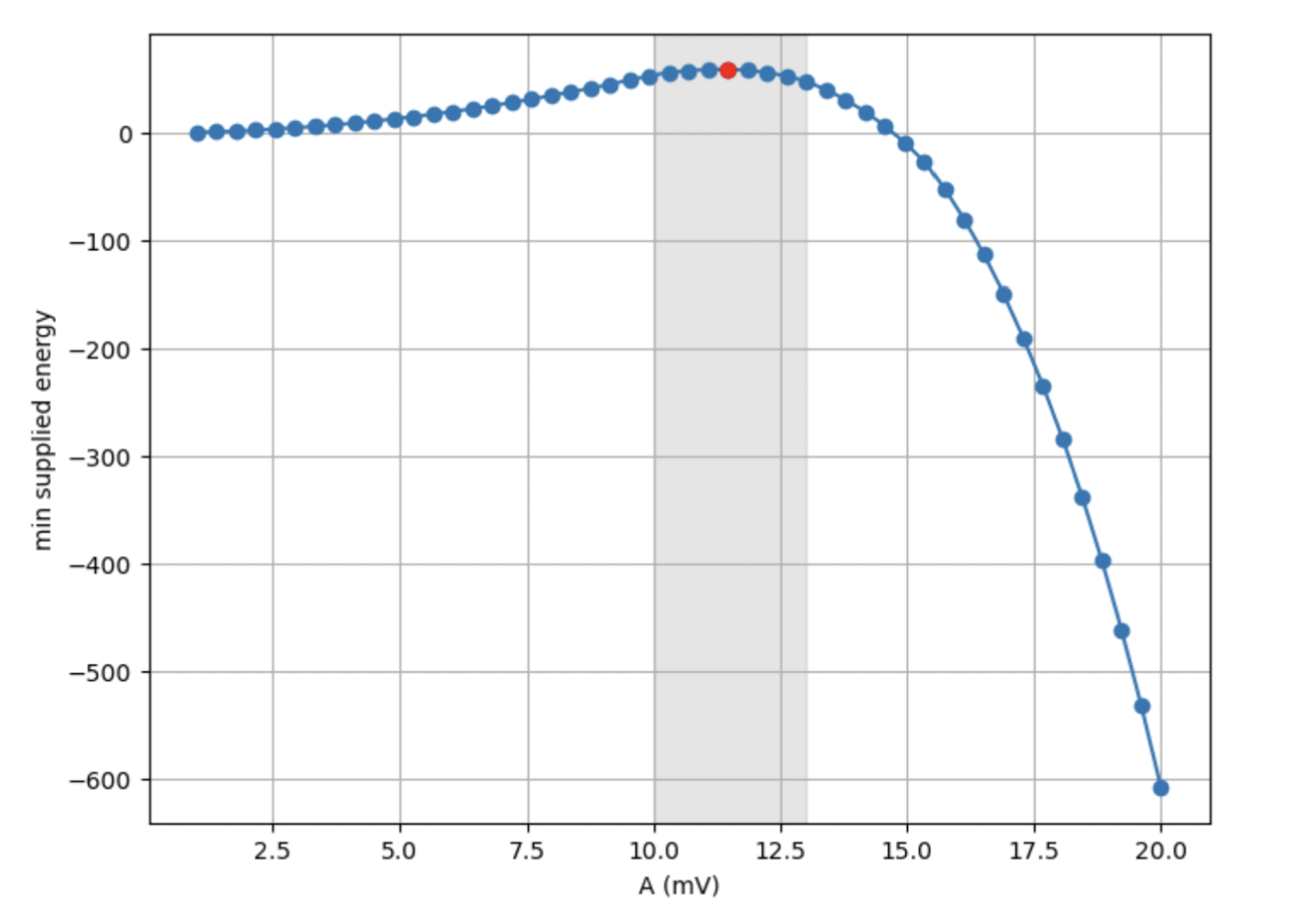} 
     \caption{The excitatory threshold of the Hodgkin--Huxley circuit.}
     \label{fig:HH-excitory threshold}
\end{figure}
\subsection{Interpretation of the excitatory threshold}
Fig.~\ref{fig:HH-excitory threshold} shows the minimal energy required as a function of the target voltage \( A \). A clear local maximum appears at around \(A\approx11.5 mV, \alpha \approx 0.62\). This local maximum is a good approximation of  the excitatory threshold of the Hodgkin--Huxley circuit.
Below the threshold voltage, approaching higher voltages requires increasingly large energy inputs due to local passive behavior of the circuit, as the system resists moving away from its stable resting equilibrium. Once the threshold is exceeded, internal energy sources facilitates a regenerative, self-sustained voltage rise (the ``spike"), leading to a spiking event without the need of additional energy supply. 

Hence, the excitatory threshold is naturally characterized as a local maximum of the required supply. This maximum represents the minimal energy barrier between subthreshold (analog-like) and suprathreshold (digital-like) behaviors in Hodgkin--Huxley model.

\subsection{Estimating the inhibitory threshold}

\begin{figure}[h]  
     \centering
     \includegraphics[width=0.47\textwidth]{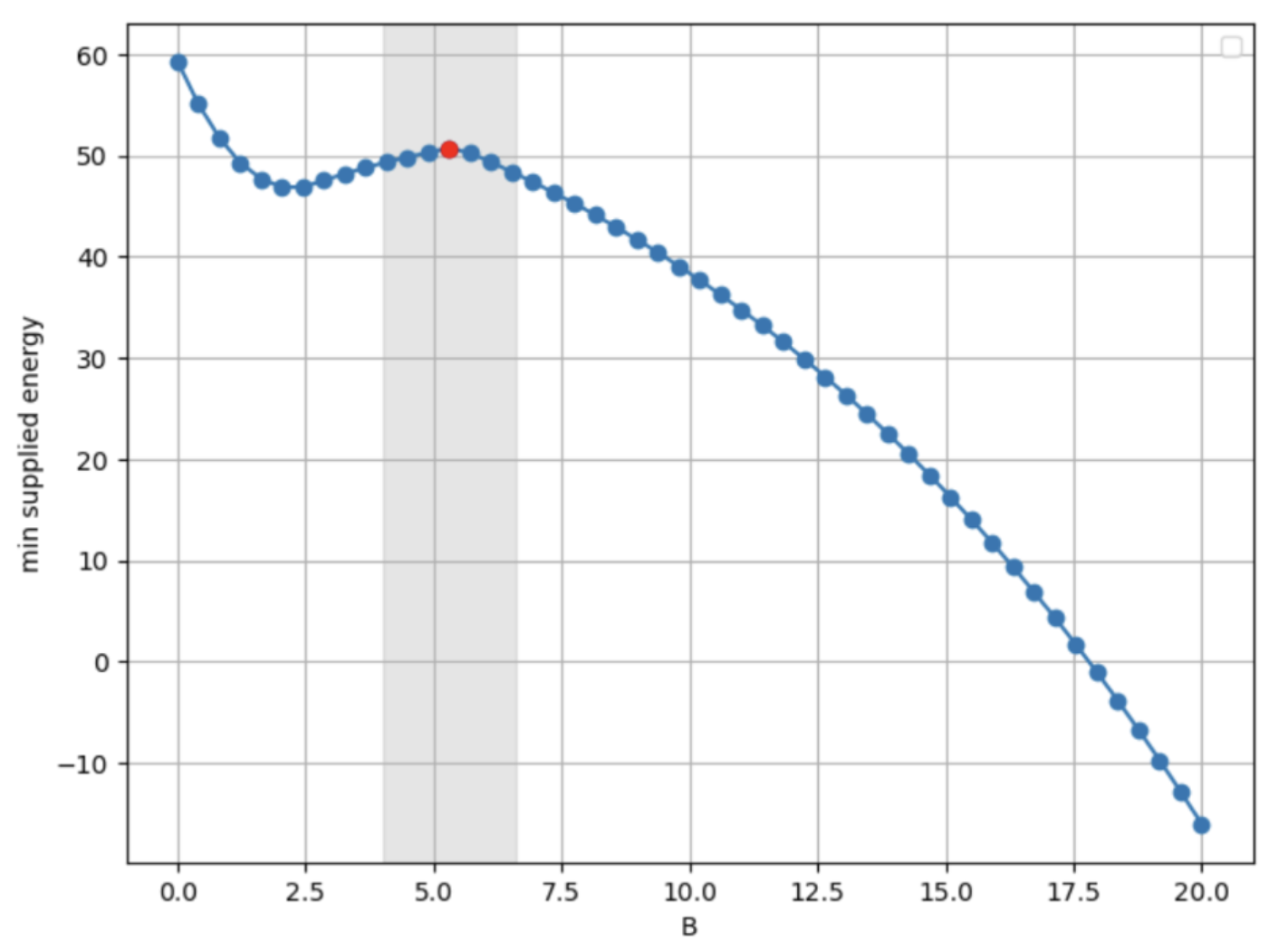} 
     \caption{Adding inhibition can lower the threshold voltage.}
     \label{fig:HH-inhibition}
\end{figure}

Fig.~\ref{fig:HH-inhibition} illustrates two distinct thresholds in Hodgkin--Huxley model: an excitatory threshold, reached with a depolarizing current, and an inhibitory threshold, reached with a hyperpolarizing current. The inhibitory threshold is responsible for the mechanism of {\it rebound} excitability, which is often neglected in neuromorphic engineering but is nevertheless a key excitability function in many neuronal behaviors. 

To approximate this second threshold using our energy-based definition, we enlarge the subspace of past voltages to allow for both excitation and inhibition:  \(Ae^{\alpha t}-Be^{\beta t}\).
The local maximum in Fig.~\ref{fig:HH-inhibition} is now computed by fixing the search space to be \(A^*e^{\alpha^* t}-Be^{\beta t}\), where \(A^*, \alpha^*\) corresponds to the local maximum in Fig.~\ref{fig:HH-excitory threshold}. The resulting threshold voltage is lower comparing to that of no inhibition.
Even though the numerical investigation above is rather limited, it provides very encouraging results. First, it captures two distinct threshold properties of the Hodgkin--Huxley model with an elementary family of past voltages. 
Second, the optimal strategy has a clear biophysical interpretation: the excitatory threshold only exploits the sodium current, which provides negative conductance at a certain time scale and amplitude. We observe a qualitative correspondence between the optimal voltage trajectory and the sodium current’s activation profile, whose activation is maximal around -50 mV and which activates with a time constant of about 0.6. Likewise, the inhibitory threshold improves on the excitatory threshold by exploiting the dissipation properties of the potassium current: by first bringing the potential close to the voltage of the potassium battery, the dissipation of the potassium current is strongly reduced, allowing to trigger an event with less supplied energy that in a purely depolarizing scenario.

The fact that the optimal trajectories can be given an easy physical interpretation is very encouraging. It suggests that the proposed energy-based characterization is robust to the details of the model and that efficient numerical procedures can be developed for more general models by exploiting the physical properties of the distinct current sources.

\section{Discussion}

We have proposed a novel definition for the threshold of excitable systems. Departing from the classical characterization of excitability by means of bifurcation analysis, our approach is instead energy-based: we define the threshold by quantifying the minimal external energy required to drive the system from its resting state to an event-triggering condition.

The energy-based characterization is expressed in the classical language of dissipativity, naturally leading to its formulation as an optimal control problem. We have demonstrated that this optimal control problem possesses a rich underlying structure, which will be further explored in a forthcoming publication that builds on recent advances in gradient modeling of memristive systems.

The illustration of the proposed energy threshold in the simplest examples of excitable models is encouraging, as it aligns with more elementary notions of thresholds---for instance, the saddle point in a bistable circuit. The numerical investigation on Hodgkin Huxley model is also encouraging, illustrating how solutions to the optimal control problem connect directly to the physical properties of the circuit elements. While preliminary, the results in this paper suggest the potential of an energy-based perspective for a general theory of excitability, one that combines physical interpretation and computational tractability.




\bibliographystyle{IEEEtran} 
\bibliography{cdc25-threshold.bib} 






\end{document}